\documentclass[12pt,reqno]{amsart}
\usepackage{fullpage}
\usepackage{amsmath,amssymb,amsthm}
\usepackage{bm}
\usepackage{times}
\usepackage{verbatim}

\newtheorem*{theorem}{Theorem}
\newtheorem{lemma}{Lemma}
\newtheorem*{proposition}{Proposition}
\newtheorem*{sconjecture}{Szpiro Conjecture}
\newtheorem*{msconjecture}{Modified Szpiro Conjecture}
\newtheorem*{uniconjecture}{Uniform $\bm{abc}$ Conjecture}
\newtheorem{abconjecture}{$\bm{abc}$ Conjecture, Version}
\newtheorem*{nconjecture}{$\bm{abc}$ Conjecture, $n$-Variable Version}
\theoremstyle{definition}
\newtheorem*{definition}{Definition}
\newtheorem{example}{Example}
\newtheorem*{remark}{Remark}

\renewcommand{\mod}[1]{{\ifmmode\text{\rm\ (mod~$#1$)}\else\discretionary{}{}{\hbox{ }}\rm(mod~$#1$)\fi}}

\DeclareMathOperator{\lcm}{lcm}
\DeclareMathOperator{\Nm}{Norm}

\newcommand{\C}{{\mathbb C}}

\newcommand{\R}{{\mathbb R}}
\newcommand{\Q}{{\mathbb Q}}

\newcommand{\Z}{{\mathbb Z}}

\newcommand{\ep}{\varepsilon}

\newcommand{\thm}{\begin{theorem}}
\newcommand{\thme}{\end{theorem}}
\newcommand{\pf}{\begin{proof}}
\newcommand{\pfe}{\end{proof}}
\newcommand{\lm}{\begin{lemma}}
\newcommand{\lme}{\end{lemma}}
\newcommand{\prop}{\begin{proposition}}
\newcommand{\prope}{\end{proposition}}
\newcommand{\cor}{\begin{corollary}}
\newcommand{\core}{\end{corollary}}
\newcommand{\df}{\begin{definition}}
\newcommand{\dfe}{\end{definition}}
\newcommand{\ex}{\begin{example}}
\newcommand{\exe}{\end{example}}
\newcommand{\rmk}{\begin{remark}}
\newcommand{\rmke}{\end{remark}}
\newcommand{\cj}{\begin{conjecture}}
\newcommand{\cje}{\end{conjecture}}
\newcommand{\abcj}{\begin{abconjecture}}
\newcommand{\abcje}{\end{abconjecture}}
\newcommand{\scj}{\begin{sconjecture}}
\newcommand{\scje}{\end{sconjecture}}
\newcommand{\mscj}{\begin{msconjecture}}
\newcommand{\mscje}{\end{msconjecture}}
\newcommand{\ucj}{\begin{uniconjecture}}
\newcommand{\ucje}{\end{uniconjecture}}
\newcommand{\ncj}{\begin{nconjecture}}
\newcommand{\ncje}{\end{nconjecture}}

\begin{document}

\begin{abstract}
The $abc$ conjecture, one of the most famous open problems in number theory, claims that three positive integers satisfying $a+b=c$ cannot simultaneously have significant repetition among their prime factors; in particular, the product of the distinct primes dividing the three integers should never be much less than~$c$. Triples of numbers satisfying $a+b=c$ are called {\em $abc$ triples} if the product of their distinct prime divisors is strictly less than~$c$. We catalog what is known about $abc$ triples, both numerical examples found through computation and infinite familes of examples established theoretically. In addition, we collect motivations and heuristics supporting the $abc$ conjecture, as well as some of its refinements and generalizations, and we describe the state-of-the-art progress towards establishing the conjecture.
\end{abstract}

\title{\boldmath $abc$ triples}
\author[Greg Martin and Winnie Miao]{Greg Martin and Winnie Miao}
\maketitle

\section{Introduction}
A, B, C \dots\ only in mathematics could such a trite trio of letters signify a major outstanding open problem with significant connections to multiple topics. The \textit{$abc$ conjecture} is a simple-to-state yet challenging problem in number theory that has stumped mathematicians for the past 30 years. It has become known for its large number of profound implications in number theory and particularly in Diophantine equations; among these myriad consequences are Fermat's last theorem (up to finitely many counterexamples), Mordell's conjecture \cite{Elk}, and Roth's theorem \cite{Bom} (see \cite{Ntj} for a more comprehensive list). The $abc$ conjecture is deeply intriguing because it unveils some delicate tension between the additive and multiplicative properties of integers, the bread and butter of number theorists.

The purpose of this article is to discuss examples and constructions of $abc$ triples, which are trios of integers demonstrating that the $abc$ conjecture, if true, must be only barely true. To do so we must first, of course, describe the $abc$ conjecture itself. We begin with a preliminary definition: the {\em radical} of an integer $n$, denoted by $R(n)$, is the product of all the distinct prime factors of $n$. For example, $600=2^4\cdot3\cdot5^2$ and so $R(600) = 2\cdot3\cdot5 = 30$. In other words, $R(n)$ is the largest squarefree divisor of $n$. The radical is a multiplicative function: in particular, for pairwise relatively prime integers $a$, $b$ and $c$, we have $R(abc) = R(a)R(b)R(c)$. We may now state (the first version of) the $abc$ conjecture, which postulates that the radical of three additively-related numbers cannot often be much smaller than the numbers themselves.

\abcj \label{finitely many c^1-ep}
For every $\ep>0$, there exist only finitely many triples $(a,b,c)$ of relatively prime positive integers satisfying $a+b=c$ for which \[R(abc)<c^{1-\ep}.\]
\abcje

A typical integer's radical is not too much smaller than the integer itself, and so $R(abc)$ is often about as large as $abc$---that is, much larger than $c$. Yet there are rare occurrences of triples $(a,b,c)$ satisfying the hypotheses of the $abc$ conjecture where $c$ is in fact greater than $R(abc)$. These special cases are referred to as {\em $abc$ triples}; the smallest such example is $(a,b,c)=(1,8,9)$, for which $R(abc)=R(36)=6<9$.

Furthermore, one can even construct an infinite sequence of $abc$ triples! One such example is $(a,b,c)=(1,9^n-1, 9^n)$: since $9^n-1\equiv1^n-1\equiv0\mod8$, we see that $8$ divides $9^n-1$ for every positive integer $n$. Writing $b=2^3 k$ for some positive integer $k$, we calculate that $R(abc)=R(a)R(b)R(c)=1\cdot R(2^3 k)\cdot 3$ is at most $2k \cdot 3 = 6k$, which is less than $c=8k+1$ for every $n$. We call this an {\em infinite family of $abc$ triples}; we will see many more infinite families in Section~\ref{sec: inf families}.

As is often the case, the literature contains various equivalent formulations of the $abc$ conjecture, a few of which we list now (others will appear as we proceed through the paper). For one thing, the $abc$ conjecture is just as commonly stated with the epsilon on the opposite side:

\abcj \label{finitely many R^1+ep}
For every $\ep>0$, there exist only finitely many triples $(a,b,c)$ of relatively prime positive integers satisfying $a+b=c$ for which \[c>R(abc)^{1+\ep}.\]
\abcje

Version \ref{finitely many c^1-ep} and Version \ref{finitely many R^1+ep} can be effortlessly obtained from each other, although we need to remember that both statements are ``for every $\ep>0$'' statements: for example, the inequality in Version~\ref{finitely many c^1-ep} with a given $\ep$ implies the inequality in Version~\ref{finitely many R^1+ep} with $\ep$ replaced by $\frac{\ep}{1-\ep}$. Different versions are more or less useful in different contexts; Version~\ref{finitely many R^1+ep}, for instance, is closely connected to the ``quality'' of an $abc$ triple, a quantity we will define in Section \ref{sec: exh comput}.

For a given $\ep$, if there are only finitely many $abc$ triples for which $R(abc) < c^{1-\ep}$, then there are only finitely many values of ${R(abc)}/{c^{1-\ep}}$ that are less than $1$, and we can choose the minimum such value and call it $K(\ep)$, say. Therefore Version~\ref{finitely many c^1-ep} of the $abc$ conjecture implies a new version:

\abcj \label{constant c^1-ep}
For every $\ep>0$ there exists a positive constant $K(\ep)$ such that all triples $(a,b,c)$ of relatively prime positive integers with $a+b=c$ satisfy \[R(abc) \geq K(\ep) c^{1-\ep}.\]
\abcje

\noindent This new formulation really is equivalent to Version~\ref{finitely many c^1-ep}---more precisely, Version~\ref{constant c^1-ep} with a given positive $\ep$ implies Version~\ref{finitely many c^1-ep} for any larger $\ep$. There is a parallel reformulation from Version~\ref{finitely many R^1+ep}:

\abcj \label{constant R^1+ep}
For every $\ep>0$ there exists a positive constant $K'(\ep)$ such that all triples $(a,b,c)$ of relatively prime positive integers with $a+b=c$ satisfy \[c \le K'(\ep) R(abc)^{1+\ep}.\]
\abcje

It might be nice to be able to leave out the hypothesis that the three integers $(a,b,c)$ are relatively prime; however, this condition is in fact indispensable. (It is worthwhile to point out the slight difference between a set of integers being \textit{relatively prime} and being \textit{pairwise relatively prime}: relatively prime means there is no common prime factor shared by all its elements, while pairwise relatively prime means that any two chosen integers from the set have no common factor. For example, the set $\{6,10,21\}$ is relatively prime but not pairwise relatively prime. Fortunately in our case, the $abc$ conjecture deals only with trios of integers related by the equation $a+b=c$; as it turns out, this equation ensures that any relatively prime set $(a,b,c)$ must also be pairwise relatively prime.) Without that hypothesis, nothing would stop us from multiplying any given triple by a huge power of a prime $p$, which would increase $c$ as much as we wanted while only increasing the radical $R(abc)$ by a factor of $p$ at most. The most extreme example of this undesirable inflation is the triple $(a,b,c) = (2^n,2^n,2^{n+1})$, for which $c=2^{n+1}$ can be made as much larger than $R(abc)=2$ as we wish.

Likewise, the epsilon appearing in the statements of the conjecture might seem like a nuisance, but it turns out to be a necessity. We have already shown to be false the more simplistic assertion that $c$ can be greater than the radical $R(abc)$ only for finitely many triples; it is even false that the ratio $c/R(abc)$ is bounded above. Section~\ref{sec: inf families} is devoted to recording several examples that refute these epsilon-less statements; many of these examples are ``well known to the experts'' yet decidedly hard to find in the literature, and we hope gathering them together here (along with citations, where known) is a helpful service to those studying this topic.

Before we take on that task, however, we spend some time in Section~\ref{sec: exh comput} looking at some numerical examples of $abc$ triples that have been garnered over the years and by examining various computational techniques of obtaining such triples. After presenting the aforementioned infinite families of $abc$ triples in Section~\ref{sec: inf families}, we then delve into the motivation behind this deep conjecture in Section~\ref{sec: backgr, motiv}. Lastly, in Section~\ref{sec: generalizations, results} we present some refinements and generalizations of the $abc$ conjecture, and discuss progress towards the conjecture and its current status. Although some of these later results and extensions are a bit technical, the large majority of the material we present is pleasantly elementary and accessible.

\section{Numerical examples of $\bm{abc}$ triples} \label{sec: exh comput}

Because the $abc$ conjecture has become so prominent in the last thirty years, corresponding roughly to the era of widespread and easily accessible computation, it is no surprise that people have developed a sustained interest in compiling numerical examples of $abc$ triples. As a matter of fact, one can go to an online $abc$ triples database~\cite{Hur} and list all $abc$ triples of positive integers up to any bound less than $10^8$, or input any integer in that range to search for $abc$ triples containing it. For instance, there are exactly seven $abc$ triples with $c=10^8$: the one with the largest value of $b$ is $(a,b,c) = (351$,$297, 99$,$648$,$703, 100$,$000$,$000) = (3^4\cdot 4$,$337,7^7\cdot11^2,2^8\cdot5^8)$, for which $R(abc) = 10$,$018$,$470$.

In fact, computations of $abc$ triples have been carried out for much larger ranges. Typically such computations record the triples they find according to their ``quality'':

\begin{definition}
Given a triple $(a,b,c)$ of relatively prime positive integers such that $a+b=c$, the {\em quality} $q(a,b,c)$ of the triple is defined to be
\[
q(a,b,c) = \frac{\log c}{\log R(abc)}.
\]
\end{definition}

\noindent For example, the quality of the smallest $abc$ triple is $q(1,8,9) = \frac{\log 9}{\log 6} = 1.22629\dots$. By this definition, a triple will be an $abc$ triple only if $q(a,b,c)>1$. And indeed, we can reformulate the $abc$ conjecture yet again, by solving the inequality in Version~\ref{finitely many R^1+ep} for $1+\ep$:

\abcj \label{quality version}
For every $\ep>0$, there exist only finitely many triples $(a,b,c)$ of relatively prime positive integers satisfying $a+b=c$ for which $q(a,b,c) > 1+\ep$.
\abcje

By looking at de Smit's website~\cite{Smt}, we see for example that among numbers with at most twenty digits, there are exactly $236$ $abc$ triples of quality at least $1.4$. Atop that list is the triple
\begin{equation} \label{reyssat}
\text{$(a,b,c) = (2,6$,$436$,$341,6$,$436$,$343)=(2, 3^{10}\cdot109, 23^5)$},
\end{equation}
for which $q(a,b,c) = 1.62991\dots$; this is the highest quality of any known $abc$ triple (and possibly the highest quality of any $abc$ triple in the universe!). This triple was discovered in 1987 by the French mathematician E. Reyssat (apparently by ``brute force'', according to \cite[page 137]{LZ}). de Smit's list also includes an $abc$ triple, discovered by I. J. Calvo, where $c$ has a whopping $2$,$958$ digits: the triple
\begin{multline}
(a,b,c) = \big( 3^3 \cdot 31^3 \cdot A, \\
5^{362}\cdot 7^{109}\cdot 11^{7}\cdot 17^{326}\cdot 37^{11}\cdot 53^{33}\cdot 59^{179}\cdot 67^{137}\cdot 79^{76}\cdot 103^{348}\cdot 109^{12}\cdot 113^{103}\cdot 131^{42}\cdot 151^{12}\cdot 163^{166} , \\
2^{465}\cdot 13^{76}\cdot 19^{57}\cdot 23^{611}\cdot 29^{19}\cdot 41^{11}\cdot 43^{98}\cdot 61^{84}\cdot 71^{13}\cdot 73^{250}\cdot 83^{30}\cdot 89^{10}\cdot 97^{80} \cdot \\
\cdot 101^{45}\cdot 127^{7}\cdot 137^{8}\cdot 139^{3}\cdot 167^{253}\cdot 173^{25} \big)
\label{big triple}
\end{multline}
has quality at least $1.01522\dots$, where $A = (c-b)/93^3$ is a number with $2$,$854$ digits. (Interestingly, as is often the case with large numbers, $A$ is easily shown to be composite---by calculating that $2^{A-1}\not\equiv1\mod A$ and invoking Fermat's little theorem, for example---but its factorization is unknown.)

{\em Reken mee met ABC}, hosted by the Mathematical Institute of Leiden University~\cite{Rkm}, is a distributive computing program aiming to collect experimental data on the $abc$ conjecture. The project is based on the BOINC platform~\cite{BOINC}, and any individual with a computer can download the software and join in the hunt for $abc$ triples. The project currently has over $150$,$000$ users and has tested nearly three quintillion triples---not too much less than the number of insects on Earth!

\vskip12pt
People have developed many different techniques for finding $abc$ triples, using tools from all parts of number theory and neighboring fields. To give a flavor of the wide variety of techniques, we describe six of them now.

\subsection{ABC@home algorithm}

The ABC@home project, which supports the {\em Reken mee met ABC} distributed computation described above, uses the following algorithm~\cite{ABC@h} to search exhaustively for $abc$ triples.

Suppose that $(a,b,c)$ is an $abc$ triple of numbers all less than~$N$. Rename the integers $\{a,b,c\}$ as $\{x,y,z\}$ so that $x$, $y$, and $z$ have the smallest, middle, and largest radical, respectively. Since $(a,b,c)$ is an $abc$ triple, we have $R(a)R(b)R(c) < c < N$, and so $R(x)R(y)R(z) < N$. From this inequality and the inequalities $R(x)<R(y)<R(z)$, it is easy to deduce that $R(y) < \sqrt{N}$ and $R(x) < N/R(y)^2$.

We may therefore search for $abc$ triples up to $N$ by sorting them according to their smallest two radicals $R(x)$ and $R(y)$, both of which are at most $\sqrt N$. First, we make a list of all of the squarefree numbers less than $\sqrt N$ (by a variant of the sieve of Eratosthenes, say). For every pair of relatively prime numbers $(r,s)$ from this list that satisfy $r<N/s^2$, we calculate all pairs of numbers $(x,y)$ for which $R(x)=r$ and $R(y)=s$. There are two ways of completing the pair $(x,y)$ to a triple where two numbers sum to the third: we can set either $z=x+y$ or $z=|x-y|$. If $s<R(z)<N/rs$, then we have discovered a new $abc$ triple, namely the sorted ordering of $(x,y,z)$.

\subsection{Continued fractions}

The {\em (simple) continued fraction} of an irrational number $\theta$ is an expression of the form
\begin{equation} \label{continued fraction}
\theta = a_0 + \frac{1}{a_1 + \frac{1}{a_2 + \frac{1}{a_3 + \cdots}}},
\end{equation}
where $a_0$ is an integer and $a_j$ is a positive integer for each $j\ge1$. Calculating the ``partial quotients'' $a_0,a_1,\dots$ of a given irrational number $\theta$ turns out to be a simple variant of the Euclidean algorithm (which at its heart is simply division with remainder). If we replace the infinite tail ${a_j + \frac{1}{a_{j+1} + \cdots}}$ of the continued fraction with $a_j$ itself, we obtain a rational number called the {\em $j$th convergent} to~$\theta$. The theory of these convergents, and how they comprise the best rational approximations to $\theta$ in a suitable sense, is extremely interesting~\cite[Chapter 7]{NZM}.

For example, we calculate the continued fraction of the irrational number $\sqrt[5]{109}$, which has been cunningly chosen for its relationship to Reyssat's example~\eqref{reyssat}:
\[\sqrt[5]{109} = 2+ \frac{1}{1+ \frac{1}{1+ \frac{1}{4+ \frac{1}{77\text{,}733 + \frac{1}{2+\cdots}}}}}. \]
Noting that the quantity ${4+ \frac{1}{77\text{,}733 + \frac{1}{2+\cdots}}}$ is extremely close to $4$, we form the approximation
\[\sqrt[5]{109} \approx  2+ \frac{1}{1+ \frac{1}{1+ \frac{1}{4}}} = \frac{23}{9},\] 
which is the third convergent to $\sqrt[5]{109}$. (In this particular case, we might have found this approximation just by examining the decimal expansion $\sqrt[5]{109} = 2.555555399...$!) This approximation tells us that $9^5  \cdot109 \approx 23^5$, and in fact their difference is exactly $2$, yielding Reyssat's triple $(2, 9^5 \cdot109, 23^5)$.

In general, we begin with an irrational root $\theta = \sqrt[k]{m}$ of an integer $m$ and compute its continued fraction. At any point, when we see a relatively large partial quotient $a_{j+1}$, we truncate the infinite continued fraction~\eqref{continued fraction} after $a_j$ to obtain the $j$th convergent, which we write as $p/q$. We have thus found integers $p$ and $q$ such that $p/q \approx \sqrt[k]m$, or equivalently $mq^k \approx p^k$. We then check the triple candidate $(|mq^k - p^k|, mq^k, p^k)$ to see whether its quality exceeds~$1$.

For the curious reader, \cite{BB} contains a list of ninety $abc$ triples, all with quality exceeding $1.4$, that can be found via this continued fraction method.

\subsection{The LLL method}
Another interesting method to find $abc$ triples, proposed by Dokchitser~\cite{Dok}, employs a famous ``lattice basis reduction'' algorithm by Lenstra, Lenstra, and Lov\'asz~\cite{LLL}. A {\em lattice} is a discrete subgroup of $\R^n$ that is closed under addition; for example, the usual integer lattice $\Z^3$ is the set of all integer linear combinations of the vectors $(1,0,0),\, (0,1,0),\, (0,0,1)$ inside~$\R^3$. Those three vectors form a {\em basis} for the integer lattice, but so do say $(12, 34, 39),\, (20, 57, 65),\,(95, 269, 309)$; just like vector spaces, a lattice can have many basis. Given a complicated basis for a lattice, like this latter one, the {\em LLL algorithm} converts it into a much nicer basis, like the former one---one with smaller entries, and for which the basis elements are nearly orthogonal.

To apply this tool to the construction of $abc$ triples, we select large integers $r,s,t$ that are comparable in size and have very small radicals (high powers of small primes, for example, or products of these). If we can find small integers $u,v,w$ such that
\begin{equation} \label{eqn: LLL}
ur+vs+wt=0,
\end{equation}
then $(|u|r,|v|s,|w|t)$ has a good chance of being an $abc$ triple: the radicals of $r,s,t$ were all chosen to be small, and the integers $|u|,|v|,|w|$ themselves are small and can only contribute so much to the radical of the product.

The set of all integer vectors $(u,v,w)$ satisfying equation~\eqref{eqn: LLL} is a two-dimensional sublattice of $\Z^3$; however, the usual methods of finding a basis for this sublattice result in basis vectors with very large entries. We run the LLL algorithm on this basis to find a reduced basis $\{\mathbf b_1, \mathbf b_2\}$ for the lattice of solutions to equation~\eqref{eqn: LLL}, where the new basis vectors have much smaller entries. We may now consider any linear combination $(u,v,w) = s_1\mathbf b_1 + s_2 \mathbf b_2$, where $s_1,s_2$ are small integers, and test the triple $(|u|r, |v|s, |w|t)$ to see if it is an $abc$ triple.

In this fashion, Dokchitser was able to obtain 41 new $abc$ triples, including $(13^{10}\cdot37^2, 3^7 \cdot19^5 \cdot71^4 \cdot223, 2^{26}\cdot5^{12}\cdot1\text{,}873)$ which has a quality of $1.5094$, the $11$th highest quality known.

\subsection{Transfer method} \label{transfer section}

Yet another approach to finding new $abc$ triples is to take existing triples and ``transfer" them, using certain polynomial identities, to create new triples.

For example, note that if $a+b=c$, then $a^2+c(b-a)=b^2$, since $c(b-a)=(b+a)(b-a)=b^2-a^2$. Note also that if $R(abc)<c$, then
\begin{equation} \label{transfer 1}
R(a^2\cdot c(b-a)\cdot b^2) \le R(a)R(b)R(c)R(b-a) = \frac{R(abc)}c c R(b-a)< c(b-a) < b^2
\end{equation}
as well. In other words, if $(a,b,c)$ is an $abc$ triple with $a<b$, then $(a^2,c(b-a),b^2)$ is also an $abc$ triple. Indeed, if the quality $q(a,b,c)$ is larger than $1$, then a quick calculation~\cite[page 16]{Hrs} shows that
\[
q(a^2,c(b-a),b^2) > \frac{2q(a,b,c)}{q(a,b,c)+1} > 1.
\]
For future reference, we also note a slight improvement: if $(a,b,c)$ is an $abc$ triple where $a$ and $b$ are both odd (which forces both $c$ and $b-a$ to be even), then
\begin{multline} \label{transfer 2}
R(a^2\cdot c(b-a)\cdot b^2) < \frac c{R(abc)} R(a^2\cdot c(b-a)\cdot b^2) \\
\le \frac c{R(abc)} R(a)R(b)R(c)R\bigg(\frac{b-a}2\bigg) < c\bigg(\frac{b-a}2\bigg) < \frac{b^2}2.
\end{multline}

When we are looking for good numerical examples, moreover, we can try this transfer method on many known $abc$ triples and hope for some extra repeated factors in $b-a$. For example, we can start with the small $abc$-triple $(7,243,250)$, whose radical is $210$ and whose quality is about $1.03261$. Using the above transfer identity leads to the triple $(7^2,250(243-7),243^2) = (49, 59\text{,}000, 59\text{,}049)$. We know from the bound~\eqref{transfer 1} that the radical of this new triple is at most $210\cdot(243-7)$. However, $243-7=2^2\cdot59$, and the factors of $2$ are dropped from the radical since $250$ is already even. Consequently, the radical of this new triple is only $210\cdot59 = 12$,$390$, and the quality of $(49, 59\text{,}000, 59\text{,}049)$ is about $1.16568$, which is quite a bit better than the original triple.

The transfer method, then, is to start with existing $abc$ triples, apply a polynomial identity to obtain a new triple, and then check for fortunate coincidences that make the new triple even better than we already knew it would be. It is an experimentation game, where different starting triples can yield results from mediocre to extremely good. In fact, we can experiment not only with the starting triple but with the polynomial identity as well! Some other examples of such polynomial transfers, which are all easily seen to be valid when $c=a+b$, include:
\begin{align*}
(b-a)^2+4ab &= c^2\\
a^3+b^3&=c(b^2-ab+a^2)\\
a^2(a+3b) + b^2(3a+b) &= c^3\\
a^3(a+2b) + c^3(b-a) &= b^3(2a+b)\\
27c^5(b-a) + a^3(3a+5b)^2(3a+2b) &= b^3(5a+3b)^2(2a+3b).
\end{align*}
Moreover, there is even a whole family of such identities
\begin{equation*}
a^{n-k} \bigg( \sum_{j=0}^k {n \choose j} a^{k-j}b^j \bigg) + b^{k+1} \bigg( \sum_{j=0}^{n-k-1} {n \choose j} a^j b^{n-k-1-j} \bigg) = c^n
\end{equation*}
which comes from splitting the binomial formula for $(a+b)^n$ at some term with index $1\le k\le n-1$. (Note that the third identity on the above list is the $n=3$, $k=1$ case of this general family.)

The interested reader can refer to \cite[Section 2.3]{Hrs} for a detailed examination of these polynomial transfers as a way of generating triples.

\subsection{An elliptic curve method}

Before describing the next method of finding examples of $abc$ triples, which was developed by van der Horst \cite{Hrs}, we say a few words about {\em elliptic curves}. For our purposes, an elliptic curve can be defined as the set of solutions of a suitable cubic equation in two variables, such as~\eqref{eqn: elliptic curve} or~\eqref{short}. That set of solutions depends, of course, on what domain we select for the variables; it turns out to be fruitful to consider the same equation with different domains, as we will see below. Certainly, elliptic curves are very fascinating in their own right (see~\cite{Was} or~\cite{Silv}, for example, where one can find all the facts about elliptic curves that we describe in this paper). For now, we need only to talk about the group structure of an elliptic curve; we will mention $j$-invariants in the next section and other elliptic curve invariants in Section~\ref{ECs}.
 
Amazingly, the points on an elliptic curve can be turned into an abelian group (once a ``point at infinity'', representing the group identity, is included) using a suitable definition of addition: three points on the elliptic curve sum to the identity precisely when they are collinear. When the variables are allowed to be  complex numbers, the resulting abelian group is isomorphic to a (two-dimensional) torus. On the other hand, if the coefficients and the variables of the cubic equation are restricted to rational numbers, then the resulting abelian group is finitely generated (this is the Mordell--Weil theorem), thus having a free part $\Z^{\text{rank}}$ and a well-understood torsion subgroup. (The rank, on the other hand, is not well understood in general, which is why it is one of the subjects of the Birch and Swinnerton--Dyer Conjecture, one of the seven Clay Mathematics Institute's Millennium Problems~\cite{BSD}).

We now describe a slight variant of van der Horst's method of searching for $abc$ triples. For any fixed integers $x_0<y_0$, set $k=y_0^3-x_0^3$ and consider the elliptic curve given by the equation
\begin{equation} \label{eqn: elliptic curve}
y^3=x^3+k,
\end{equation}
where the variables $x$ and $y$ are allowed to be not just integers but rational numbers in general. Whenever $(x, y)=(\frac pd, \frac qd)$ is a point on this elliptic curve (for simplicity we assume that $p$, $q$, and $d$ are positive), we have $q^3=p^3+kd^3$. Clearly $R(p^3,kd^3,q^3) \le kdpq < kdq^2$, and so this triple is an $abc$ triple whenever $q>kd$, or equivalently when $y>k$; indeed, the larger $y$ is, the higher the quality of the triple will be.

It is probably not the case that $y_0$ itself is larger than $k$; however, we can use the group operation on the elliptic curve to search for rational solutions to equation~\eqref{eqn: elliptic curve} other than $(x_0,y_0)$. Simply adding the point $(x_0,y_0)$ to itself repeatedly (adding, that is, using the group law on the elliptic curve) yields a sequence of points on the elliptic curve that is typically infinite. van der Horst even develops a way of predicting which elements of this sequence will have large $y$-values: he writes down a group homomorphism from the elliptic curve to the unit circle in the complex plane that takes points with large coordinates to complex numbers near~$1$. Since it is easy to calculate which powers of a complex number are close to $1$, one can take the corresponding multiples of $(x_0,y_0)$ back on the elliptic curve and check how good the corresponding triples' qualities are. One feature of this method is that all three numbers in the $abc$ triples it generates have small radicals, not just one or two of them.

The exact algorithm and variants used by van der Horst~\cite[Sections 4.2--4.3]{Hrs} discovered some notable $abc$ triples. The point $(x,y) = (\frac{19}{93}, \frac{289}{93})$ on the elliptic curve $y^3 = x^3+30$ does not have $y>30$, but fortunately the numerator of $y$ happens to be a square, and so we get to divide the radical by an extra factor of~$17$. The resulting $abc$ triple $(19^3, 30\cdot93^3, 289^3) = (6\text{,}859, 24\text{,}130\text{,}710, 24\text{,}137\text{,}569)$ has radical $300$,$390$ and quality about $1.34778$. Moreover, the algorithm often finds rational solutions with huge numerators and denominators, and is thus suited for finding enormous $abc$ triples; van der Horst reports~\cite[Chapter 5]{Hrs} finding a point on the elliptic curve $y^3 = x^3+854$ that yields an $abc$ triple with quality about $1.01635$, where the largest integer in the triple has $340$ digits.

\subsection{Differences of \boldmath $j$-invariants} \label{j section}

We conclude this section with some exotic $abc$ triples that are found unexpectedly when discussing factorizations of ``$j$-invariants''. 

There is a beautiful link between lattices and elliptic curves: through two ``elliptic functions'' studied by  Weierstrass, it is known that every elliptic curve can be represented as $y^2=4x^3-g_2(\tau)x-g_3(\tau)$, where $g_2(\tau)$ and $g_3(\tau)$ are invariants that correspond to a fixed lattice. More specifically, they are the modular forms
\begin{align*}
g_2(\tau) &= 60 \sum_{\substack{m,n\in\Z \\ (m,n)\ne(0,0)}} \frac{1}{(m\tau+n)^4} \\
g_3(\tau) &=140 \sum_{\substack{m,n\in\Z \\ (m,n)\ne(0,0)}} \frac{1}{(m\tau+n)^6}
\end{align*}
where $\tau$, a complex number with positive imaginary part, determines the relevant lattice as the set of all numbers of the form $m\tau +n$ with $m,n$ integral. (This lattice, by the way, is exactly the lattice one needs to quotient the complex plane by to realize the elliptic curve; since a plane modulo a lattice is a torus, this description corroborates the fact that every elliptic curve is isomorphic to a torus, as mentioned in the previous section.)

Now, we define the {\em $j$-invariant} $j(\tau)$ of an elliptic curve by the formula
\[
j(\tau) = \frac{1728g_2^3(\tau)}{g_2^3(\tau)-27g_3^2(\tau)}.
\]
This $j$-invariant is a modular function with ubiquitous remarkable properties and applications in complex analysis, algebraic number theory, transcendence theory, and so on. When the argument $\tau$ lies in an imaginary quadratic field $\Q(\sqrt{-d})$ for some positive integer $d$, the values $j(\tau)$ are called ``singular moduli'', and the associated elliptic curves possess extra endomorphisms and are said to have ``complex multiplication''. This singular modulus is an algebraic integer lying in some abelian extension of $\Q(\sqrt{-d})$; remarkably, the degree of its minimal polynomial is exactly the ``class number'' $h(-d)$, which is the number of binary quadratic forms $ax^2+bxy+cy^2$ of discriminant $-d$ that are not equivalent to one another under linear changes of variables. In particular, by the Stark--Heegner theorem~\cite[Appendix C, Section 11]{Silv}, there are only thirteen negative discriminants $-d$ that have class number equal to $1$, namely $-3$, $-4$, $-7$, $-8$, $-11$, $-12$, $-16$, $-19$, $-27$, $-28$, $-43$, $-67$, and $-163$; the corresponding $j$-invariants are thus actual integers.

As it happens, these thirteen special $j$-invariants are all forced to be perfect cubes of integers. Equally marvelously, the difference of two of these special $j$-invariants is very nearly a perfect square~\cite{GrsZag,LtrVir}. The corresponding triple of integers is therefore a prime candidate for an $abc$ triple (at least, once the three integers are divided by their greatest common divisor). Gross and Zagier~\cite{GrsZag} cite an example with $\tau = (-1+i\sqrt{163})/2$, where the three integers
\begin{align*}
\frac{j(i)}{1728} &= 1 \\
\frac{-j(\tau)}{1728} &= \text{$151$,$931$,$373$,$056$,$000$} = 2^{12}\cdot5^3\cdot23^3\cdot29^3 \\
\frac{j(i)-j(\tau)}{1728} &= \text{$151$,$931$,$373$,$056$,$001$} = 3^3\cdot7^2\cdot11^2\cdot19^2\cdot127^2\cdot163
\end{align*}
form an $abc$ triple with quality about $1.20362$. Going through all $\binom{13}2=78$ possible pairs of special $j$-invariants, we find that the best resulting $abc$ triple comes from both $j(\tau_4)-j(\tau_{43})$ and $j(\tau_{16})-j(\tau_{67})$, where $\tau(d) = \frac12(d+\sqrt{-d})$: the triple is $(1, 512$,$000, 512$,$001)=(1,2^{12}\cdot5^3,3^5\cdot 7^2\cdot 43)$  and has quality about $1.44331$.

\section{Infinite families of $\bm{abc}$ triples} \label{sec: inf families}

All of the numerical examples from Section~\ref{sec: exh comput}, however interesting, cannot shed any light on whether the $abc$ conjecture is true or false: the ``only finitely many'' or ``there exists a constant'' clauses in its various versions preclude us from drawing conclusions from any finite number of examples. For that matter, any finite number of examples cannot rule out even more ambitious possible versions of the $abc$ conjecture. For instance, could there be an absolute constant $S>0$ such that $c < S\cdot R(abc)$ always? This statement, similar to the $abc$ conjecture but without the messy epsilons, might be called the ``simplistic $abc$ conjecture''. Again, no finite amount of computation can resolve this question.

What we need, to help us decide whether these statements are true or false, are constructions of infinite families of $abc$ triples. And it turns out that several such constructions exist; any one of these constructions shows that the simplistic $abc$ conjecture is false. In other words, the constructions in this section demonstrate that the epsilons in the $abc$ conjecture are necessary if we hope that the assertion is true.

\subsection{The transfer method again}

Recall from Section~\ref{transfer section} that if $(a,b,c)$ is an $abc$ triple, then so is $(a^2 , c(b-a) , b^2)$. In particular, if $(1,c-1,c)$ is an $abc$ triple, then so is $(1,c^2-2c,(c-1)^2)$. Of course, we can iterate this transfer multiple times in a row: for example, $(1,(c^2-2c)^2-1,(c^2-2c)^2) = (1,c^4-4 c^3+4 c^2-1,c^4-4 c^3+4 c^2)$ will also be an $abc$ triple. As it happens, doing this double transfer always allow us to remove an extra factor of $2$ from the radical. For example, suppose that $c$ is odd. Then, by setting $a=1$ and $b=c-1$ in the third and last terms of the chain of inequalities~\eqref{transfer 1}, we know that
\[
\frac{(c-1)^2}{R\big((c^2-2c)(c-1)^2\big)} \ge \frac c{R\big((c-1)c\big)}.
\]
But now $(c-1)^2$ is even, so replacing $a$, $b$, and $c$ in the second and last terms of the chain of inequalities~\eqref{transfer 2} with $1$, $c^2-2c$, and $(c-1)^2$, we find that
\begin{equation} \label{double transfer}
\frac{(c^2-2c)^2}{R\big(((c^2-2c)^2-1)(c^2-2c)^2\big)} \ge 2\frac{(c-1)^2}{R\big((c^2-2c)(c-1)^2\big)} \ge 2 \frac c{R\big((c-1)c\big)}.
\end{equation}

We can iterate this double transfer endlessly to create an infinite sequence. Let us set $c_0=9$, corresponding to the $abc$ triple $(1,8,9)$, and for every $n\ge0$ define $c_{n+1} = c_n^4-4 c_n^3+4 c_n^2$. For example, $c_1=3$,$969$, corresponding to the double transfer $(1,8,9) \rightarrow (1,63,64) \rightarrow (1,3$,$968,3$,$969)$. Equation~\eqref{double transfer} tells us that
\[
\frac{c_{n+1}}{R\big((c_{n+1}-1)c_{n+1}\big)} \ge 2 \frac{c_n}{R\big((c_n-1)c_n\big)}
\]
for every $n\ge0$. Since ${c_0/R((c_0-1)c_0)} = \frac32$, this implies that
\begin{equation} \label{broke simplistic}
\frac{c_n}{R\big((c_n-1)c_n\big)} \ge 2^{n-1}\cdot3
\end{equation}
for every $n\ge0$. And since $2^{n-1}\cdot3$ exceeds any constant we might care to name in advance, we have just created an infinite sequence of $abc$ triples $(1,c_n-1,c_n)$ that repudiates the ``simplistic $abc$ conjecture''!

We can convert the inequality~\eqref{broke simplistic} into a quantitative measure of how much smaller than $c_n$ this radical is. Note that $c_n\le c_{n-1}^4$ for every $n\ge1$, and so $c_n \le c_0^{4^n} = 9^{4^n}$. In particular, $\log c_n \le 4^n \log 9$, and so $2^n \ge \sqrt{\log c_n}/\sqrt{\log 9}$. It now follows from~\eqref{broke simplistic}, when $a=1$, $b=c_n-1$, and $c=c_n$, that
\begin{equation} \label{first R bound}
R(abc) \le \frac c{2^{n-1}\cdot3} \le \frac{2\sqrt{\log9}}3 \frac c{\sqrt{\log c}}.
\end{equation}
To this point, it hasn't mattered which logarithm we've been using, but now we clarify that we are using $\log x$ to denote the natural logarithm (which is often written $\ln x$), as is standard in analytic number theory. With that admission out of the way, we remark that the constant $\frac{2\sqrt{\log9}}3$ is approximately $0.988203$.

This bound for the radical of these triples can be re-expressed as an inequality about their quality: the lower bound
\begin{align}
q(a,b,c) = \frac{\log c}{\log R(abc)} &\ge \frac{\log c}{\log c - \log\sqrt{\log c} + \log(\frac23\sqrt{\log9})} \notag \\
&\ge \frac{\log c}{\log c - \frac12\log\log c} > \frac{\log c + \frac12\log\log c}{\log c} = 1 + \frac{\log\log c}{2\log c}  \label{quality bound}
\end{align}
holds when $(a,b,c) = (1,c_n-1,c_n)$. Notice that these qualities are all greater than $1$, but the lower bound does tend to $1$ as $c$ becomes larger and larger. If the lower bound tended to a constant larger than $1$, this sequence would disprove the actual $abc$ conjecture (specifically Version~\ref{quality version}) and this whole paper would need to be rewritten!

\subsection{Folklore examples} \label{folklore sec}

There are several known constructions of infinite sequences of $abc$ triples, each of which provides a counterexample to the ``simplistic $abc$ conjecture''. We present a few of these constructions in this section. Unlike the recursive construction from the previous section, these constructions have very simple closed forms which make it obvious that the smallest and largest numbers in the triples have extremely small radicals. In each case, a quick number theory lemma is required to show that the radical of the middle number is somewhat smaller than the number itself. These constructions are simple enough (to those well-versed in the field) that it is nearly impossible to determine who first came up with them; indeed, some cannot even be found explicitly in any publication despite that they are ``well known''! Part of the motivation for this paper was to ensure that these families of $abc$ triples are explicitly recorded in the literature; we have included earlier citations whenever we could locate them.

\lm
\label{lm: Euler}
If $p$ is an odd prime, then $p^2$ divides $2^{p(p-1)}-1$.
\lme

\pf
Euler's theorem \cite[page 63]{Euler} says that if $a$ and $m$ are relatively prime positive integers, then $a^{\phi(m)} \equiv 1\mod m$, where $\phi(m)$ is the Euler phi-function. Applied with $a=2$ and $m=p^2$, for which $\phi(m)=p(p-1)$, Euler's theorem yields $2^{p(p-1)} \equiv 1\mod{p^2}$, which is exactly the conclusion of the lemma.
\pfe

The following construction was recorded by Granville and Tucker~\cite{Grv}.

\ex
\label{ex: Grv}
For any odd prime $p$, set $(a,b,c) = (1, 2^{p(p-1)}-1,2^{p(p-1)} )$. We know by Lemma~\ref{lm: Euler} that $p^2$ divides $b$, and so $R(b) \le b/p$. It follows that
\begin{equation} \label{eqn: Granville's ex}
R(abc) = R(a)R(b)R(c) \leq 1 \cdot \frac{b}{p} \cdot 2 < \frac{2c}{p}.
\end{equation}
Since the sequence of primes $p$ becomes larger than any constant we want, this family of triples does contradict the ``simplistic $abc$ conjecture''.

For easier comparison to other examples, we can rewrite the right-hand side in a form involving only~$c$. Since $c<2^{p^2}$, we have $\log c < p^2\log 2$ and so $p > \frac{\sqrt{\log c}}{\sqrt{\log 2}}$. Combining this with (\ref{eqn: Granville's ex}) yields
\[
R(abc) < 2\sqrt{\log 2} \frac{c}{\sqrt{\log c}}.
\]
This upper bound for the radical has the same shape as the bound in ~equation~\eqref{first R bound} for our first example, but with the slightly worse constant $2\sqrt{\log2} \approx 1.66511$.
\exe

Our next infinite family involves a lemma providing divisibility by high powers of a prime, rather than just its square.

\lm
\label{lm: 7^(n+1) divs}
If $n$ is a nonnegative integer, then $7^{n+1}$ divides $8^{7^n}-1$.
\lme

\pf
We proceed by induction; the base case $n=0$ is immediate. Assuming the lemma is true for a particular $n$, we write
\[
8^{7^{n+1}}-1 = 8^{7\cdotp{7^n}}-1 = (8^{7^n}-1)(8^{6\cdotp{7^n}}+8^{5\cdotp{7^n}}+ \dots + 8^{7^n}+1).
\]
On the right-hand side, the first factor is divisible by $7^{n+1}$ by the induction hypothesis, while the second factor is divisible by $7$ since each of its seven terms is congruent to $1 \mod 7$. Therefore $7^{n+1}\cdot7$ divides the left-hand side, which is the statement of the lemma for $n+1$ as required.
\pfe

\ex
\label{ex: BMT}
For any nonnegative integer $n$, set $(a,b,c) = (1, 8^{7^n}-1, 8^{7^n})$. Equipped with Lemma~\ref{lm: 7^(n+1) divs}, we deduce that $R(b) \le b/7^n$ and thus
\[
R(abc)=R(a)R(b)R(c) \le 1\cdot\frac{b}{7^n}\cdot2 < \frac{2c}{7^n}.
\]
Again we have disproved the ``simplistic $abc$ conjecture'', and again we can write the right-hand side as an expression in $c$ alone, since $\log c=7^n\log8$:
\[
R(abc)<2\log 8\frac{c}{\log c}.
\]
Note that we have improved the order of magnitude of the upper bound on the radical, from the previous examples' $c/\sqrt{\log c}$ to $c/\log c$.
\exe

Variants of this construction abound. It is equally easy to prove by induction that $2^{n+2}$ divides $3^{2^n}-1$ for any $n\ge1$, and so a similar construction (attributed in~\cite[pages 40--41]{Lang} to Jastrzebowski and Spielman) with the triple $(a,b,c)=(1, 3^{2^n}-1, 3^{2^n})$ results in the upper bound
\[
R(abc)<\frac{3c}{2^{n+1}} = \frac{3\log 3}{2} \frac{c}{\log c}.
\]
Here the leading constant $\frac{3\log3}2\approx 1.64792$ is even better than $2\log8\approx 4.15888$.

Various constructions of this type are easily found by replacing $8^{7^n}$ or $3^{2^n}$ with a sequence of the form $q^{p^n}$, where $p\ge2$ is an integer dividing $q-1$. When $p$ is a prime and $q$ is a prime power, this construction was given by Stewart~\cite[Theorem 1, (3)]{Stw}. All of these constructions show that the radical is less than some constant (depending on the parameters chosen) times $c/\log c$. Moreover, the same sort of argument that led to equation~\eqref{quality bound} shows that the qualities of the $abc$ triples arising from Example~\ref{ex: BMT} are essentially as large as $q(a,b,c) > 1 + \frac{\log\log c}{\log c}$, without the factor of $2$ in the denominator. (The same bound will hold for the rest of the examples in Section~\ref{sec: inf families}.)

\medskip

Our last example differs from the previous ones: the radical of the middle number of the triple is small because high powers of several primes divide it, not just a high power of a single prime.

\lm
\label{lm: lcm divs}
For any positive integer $n$, define $L = \lcm [1,2,\dots n]$ and $t=\lfloor\frac{\log n}{\log 2}\rfloor$, and let $P = \prod_{3\le p\le n}p$ be the product of all the odd primes up to~$n$. Then $PL/2^t$ divides $2^L-1$. In particular,
\[
R(2^L-1) \le \frac{2^t(2^L-1)}L.
\]
\lme

\pf
Given an odd prime $p \le n$, let $r = \lfloor\frac{\log n}{\log p}\rfloor$, so that $p^r$ is the largest power of $p$ not exceeding~$n$. Clearly both $p^r$ and $p-1$, being at most $n$ in size, divide $L$; since they are relatively prime, their product $p^r(p-1)$ also divides~$L$. As $2$ is relatively prime to $p^{r+1}$, Euler's theorem tells us that $2^{\phi(p^{r+1})} = 2^{p^r(p-1)} \equiv 1\mod{p^{r+1}}$, and therefore $2^L \equiv 1\mod{p^{r+1}}$ since $L$ is a multiple of $p^r(p-1)$. Therefore $p^{r+1}$ divides $2^L-1$ for every odd prime $p\le n$. All of these prime powers are relatively prime to one another, and hence their product
\[
\prod_{3\le p\le n} p^{r+1} = \prod_{3\le p\le n} p \prod_{3\le p\le n} p^r = P\frac L{2^t}
\]
also divides $2^L-1$, as claimed. In this last equality, we used the fact that $\lcm [1,2,\dots n]$ is composed exactly from the highest power of each distinct prime factor found amongst the factorizations of the numbers 1 through $n$. In other words $L = 2^{\lfloor\frac{\log n}{\log 2}\rfloor} \prod_{3\le p\le n} p^{\lfloor\frac{\log n}{\log p}\rfloor} = 2^t\prod_{3\le p\le n} p^r$.

Note also that every prime dividing $L/2^t$ is an odd prime not exceeding $n$, hence divides $P$ as well. The above argument shows that $P$ divides the quotient $(2^L-1)/(L/2^t)$, and so the primes dividing $L/2^t$ are already represented in this quotient; consequently, the radical of $2^L-1$ is no larger than $(2^L-1)/(L/2^t)$.
\pfe

\ex \label{ex: lcm}
For any positive integer $n$, define $L = \lcm [1,2,\dots n]$, and set $(a,b,c) = (1,2^L-1,2^L)$. Using the notation $t=\lfloor\frac{\log n}{\log 2}\rfloor$ and $P = \prod_{3\le p\le n}p$ from Lemma~\ref{lm: lcm divs}, we have $2^t < n$ and $\log c = L\log 2$ and thus
\[
R(abc) \le 1\cdot \frac{2^t(2^L-1)}L \cdot 2 < 2\log 2\frac{nc}{\log c}.
\]
It is a bit harder than in the previous examples to write the right-hand side solely in terms of $c$, since the relationship between $n$ and $c$ is more complicated. The Chebyshev function $\psi(n) = \log \lcm[1,\dots,n]$ (often written in terms of the ``von Mangoldt function'' $\Lambda(n)$) satisfies $\psi(n) \sim n$ by the famous prime number theorem \cite[pages 74--75]{PNT}. Therefore $\log\log c = \log L + \log\log 2 = \psi(n) + \log\log 2 \sim n$, and hence we have the asymptotic inequality
\[
R(abc) \lesssim 2\log 2\frac{c\log\log c}{\log c},
\]
which has a slightly worse order of magnitude than the last two examples. For what it's worth, we can remove a factor of $2$ from the right-hand side by restricting $n$ to be just less than a power of~$2$.
\exe

\subsection{A curious divisibility}

All of the $abc$ triples constructed in this section so far share the property that their smallest number equals~$1$. However, we have a final construction to describe, one that was discovered only recently~\cite{BMT}, which has the feature that all three numbers in the constructed $abc$ triples are nearly the same size. This construction relies on the following quite strange divisibility relationship.

\lm
\label{lm: mod 6 divs}
For any positive integer $n$ satisfying $n \equiv 2 \mod 6$,
\[
\bigg(\frac{n^2-n+1}{3}\bigg)^2 \text{\quad divides \quad} n^n-(n-1)^{n-1}.
\]
\lme

\noindent Setting $n=6k+2$ for a nonnegative integer $k$ reveals that the lemma is equivalent to the curious statement:
\begin{equation} \label{k version}
(12k^2+6k+1)^2 \text{\quad divides \quad} (6k+2)^{6k+2} - (6k+1)^{6k+1}.
\end{equation}

\pf 
Given a nonnegative integer $k$, set $Q = 12k^2+6k+1$. To establsh the divisibility~\eqref{k version}, we need to show that $(6k+2)^{6k+2} \equiv (6k+1)^{6k+1} \mod{Q^2}$. Our main tool will be the following observation: if $a \equiv bQ+1 \mod{Q^2}$, then $a^j \equiv jbQ+1 \mod{Q^2}$ for any positive integer~$j$. This observation follows from the binomial expansion
\[
a^j \equiv (1+bQ)^j = \sum_{i=0}^j \binom ji (bQ)^i \equiv 1 + j\cdot bQ + \sum_{i=2}^j 0 \mod{Q^2}.
\]

Since $(6k+1)^3 = 18kQ+1$ and $-(6k+2)^3 = -(18k+9)Q+1$, we can certainly say that
\begin{align*}
(6k+1)^3 &\equiv 18kQ+1 \mod{Q^2}\\
-(6k+2)^3 &\equiv -(18k+9)Q+1 \mod{Q^2}.
\end{align*}
Raising both congruences to the $2k$th power using our observation, we see that
\begin{align*}
(6k+1)^{6k} &\equiv 2k\cdot18kQ+1 = (3Q-(18k+3))Q+1 \equiv -(18k+3)Q+1 \mod{Q^2} \\
(6k+2)^{6k} &\equiv -2k(18k+9)Q+1 = (-3Q+3)Q+1 \equiv 3Q+1 \mod{Q^2}.
\end{align*}
We now calculate that
\begin{align*}
(6k+2)^{6k+2}-(6k+1)^{6k+1} &\equiv (3Q+1)(6k+2)^2 - (-(18k+3)Q+1)(6k+1) \\
&= (3Q+1)(3Q+6k+1) + (9Q-18k-6)Q - (6k+1) \\
&= 18Q^2 \equiv0\mod{Q^2},
\end{align*}
which is what we needed to show.
\pfe

\rmk
Although Lemma~\ref{lm: mod 6 divs} has the elementary (if unilluminating) proof just given, there is in fact a deeper explanation~\cite[Proposition 4.3]{BMT} behind this interesting divisibility. It is related to the trinomial $x^n + x + 1$, which is reducible when $n\equiv2\mod6$, and the relationship between its discriminant $n^n-(n-1)^{n-1}$ and the resultant of its irreducible factors.
\rmke

All this work allows us to establish a bound for the radical of $b$ in the infinite family of $abc$ triples we will now construct.

\ex
\label{ex: (n-1)^(n-1)}
For any odd integer $k\ge7$, set $n=2^k$ and
$$
(a,b,c)=\big((n-1)^{n-1}, n^n-(n-1)^{n-1}, n^n\big).
$$
Since $n$ is congruent to $2\mod 6$, Lemma~\ref{lm: mod 6 divs} tells us that $(\frac{n^2-n+1}{3})^2$ divides~$b$. Therefore
\begin{equation} \label{Pomerance will improve}
R(abc) = R(a)R(b)R(c) \leq (n-1) \cdot \frac{b}{({n^2-n+1)/3}} \cdot 2 < \frac{6b}{n} < \frac{6c}{n}.
\end{equation}
Seeking a lower bound on $n$, we write $\log c = n\log n$ and $\log \log c = \log n + \log \log n < \frac43\log n$ when $n>100$, hence $n = \log c/\log n > \log c/(\frac34\log\log c)$ and so
\[
R(abc) < \frac{6b}{n} < \frac{6c}{\log c/(\frac34\log\log c)} = \frac{8c\log \log c}{\log c}
\]
when $k\ge7$.
\exe

As stated so far, this construction yields a bound on the radical comparable to the bound from Example~\ref{ex: lcm}, but with a worse constant (although for large $n$, the $8$ can essentially be replaced by a~$6$). However, if we choose specific values for $n$ in the previous example in a manner suggested by Carl Pomerance, we can further decrease the radicals of the corresponding $abc$ triples to be on par with the bound from Example~\ref{ex: BMT}.

\ex
\label{ex: (n-1)^(n-1) with n=8^7^j }
For any positive integer $j$, set $k=3\cdot 2^j$ in the triple of Example~\ref{ex: (n-1)^(n-1)}, so that $n=8^{7^j}$. Using Lemma~\ref{lm: 7^(n+1) divs}, we see that $7^{j+1}$ divides $n-1$ and thus $R(a)=R(n-1) \le (n-1)/7^j$. Therefore for the $abc$ triples
\[
(a,b,c)=\Big( \big(8^{7^j}-1\big)^{8^{7^j}-1}, 8^{7^j 8^{7^j}}-\big(8^{7^j}-1\big)^{8^{7^j}-1}, 8^{7^j 8^{7^j}} \Big)
\]
we may improve the bound~\eqref{Pomerance will improve} to
\[
R(abc) = R(a)R(b)R(c) \leq \frac{n-1}{7^j} \cdot \frac{b}{({n^2-n+1)/3}} \cdot 2 < \frac{6b}{7^jn} < \frac{6c}{7^jn} = 6 \log 8\frac{c}{\log c}.
\]
\exe

We end this section with a question: we have seen several elementary constructions of infinite families of $abc$ triples, all of which yield an upper bound on $R(abc)$ somewhere between $c/\sqrt{\log c}$ and $c/\log c$ in magnitude. Is there an elementary construction of a sequence of $abc$ triples satisfying $R(abc) < c/(\log c)^\lambda$ for some $\lambda>1$, or equivalently, satisfying $q(a,b,c) \gtrsim 1+{\lambda\log\log c/\log c}$? (We will see in Section~\ref{ST sec} that such sequences exist, but the proof does not supply a formula for them, merely a proof of their existence.)

\section{Background, motivation, and support for the $\bm{abc}$ conjecture} \label{sec: backgr, motiv}

The $abc$ conjecture was proposed in 1985 by Masser and Oesterl\'{e} \cite{Msr,Ost}, who were motivated by two analogous problems concerning polynomial rings and elliptic curves. In addition, after the $abc$ conjecture's appearance, number theorists found a probabilistic heuristic that also supports its statement. In this section we describe these links between the $abc$ conjecture and other branches of mathematics.

\subsection{The Mason--Stothers theorem}

Despite their very different appearances, the integers $\Z$ and the ring of polynomials with complex coefficients $\C[x]$ have a lot in common. In both settings, all nonzero elements enjoy unique factorization into irreducible elements: every integer can be written uniquely as a product of primes (and possibly $-1$), while every polynomial can be written uniquely as a product of monic linear factors $x-\rho$ (and possibly a nonzero leading coefficient in $\C$). Indeed, each ring is a principal ideal domain (PID), which is even stronger than being a unique factorization domain (UFD). In particular, one can define the radical $R(a)$ of a polynomial $a(x)\in \C[x]$ to simply be the product of all distinct monic linear factors that divide it, in perfect analogy with the radical of an integer. Similarly, one can define the greatest common divisor of two polynomials and hence decide whether two polynomials are relatively prime. (For these definitions, we ignore the leading coefficients, which are ``units'' in $\C[x]$, just as we might take absolute values of integers to ignore their sign for the purposes of examining their factors.) It follows that the degree of the radical of a polynomial in $\C[x]$ is the same as the number of distinct complex roots of the polynomial.

The integers generate the rational numbers $\Q$, which are quotients of one integer by a second nonzero integer; the polynomials generate the aptly named rational functions $\C(x)$, which are quotients of one polynomial by a second polynomial that is not identically zero. The rational numbers form the simplest example of a {\em number field} (we will say more about number fields in Section~\ref{number field sec}), while the field of rational functions over $\C$ form a {\em function field}; and it is a robust phenomenon in number theory (see for example \cite[Chapter 1, Section 14]{Neu}) that most results in number fields have analogous formulations in function fields. We have seen that irreducible polynomials correspond to prime numbers; another entry in the ``dictionary'' between the two rings is that the degree of a polynomial corresponds to the logarithm of a positive integer.

Masser's description of the $abc$ conjecture was motivated by the following theorem in the ``function field case'', independently discovered by Stothers and Mason~\cite{Sto,Msn} in the 1980s:

\begin{theorem}[Mason--Stothers]
Let $a(x), b(x), c(x) \in \C[x]$ be relatively prime polynomials satisfying $a(x)+b(x)=c(x)$. Then
\[
\max\big\{\deg(a), \deg(b), \deg(c)\big\} \le \deg(R(abc))-1.
\]
\end{theorem}

As it happens, the proof of the Mason--Stothers theorem is actually quite elementary. Some versions of the proof (see for example~\cite[Chapter IV, Sections 3 and 9]{Lang2}) rely on one important feature of polynomials that is completely absent from the integers: the ability to take derivatives. For example, it is not hard to show that a polynomial is squarefree (that is, has no repeated factors in its factorization into linear polynomials) if and only if it is relatively prime to its derivative. Number theorists would love to be able to detect squarefree integers so easily!

What would happen if, in the Mason--Stothers theorem, we translated from the function field setting to the number field setting by replacing degree with logarithm everywhere? We would obtain the statement $\max\big\{\log(a), \log(b), \log(c)\big\}+1 \le \log(R(abc))$, which, after exponentiating, becomes $e\max\{a,b,c\} \le R(abc)$, or simply $R(abc) \ge ec$ if we order the three positive integers so that $a+b=c$. This is an instance of the ``simplistic $abc$ conjecture'' we disproved thoroughly in Section~\ref{sec: inf families}. So the analogy between function fields and number fields, while fruitful, should always be taken with an epsilon grain of salt.

\subsection{The Szpiro conjecture}  \label{ECs}

In addition to the analogy with triples of polynomials, Oesterl\'e's motivation for formulating the $abc$ conjecture had an additional source: the subject of elliptic curves. We need to give a quick crash course in invariants of elliptic curves before stating the Szpiro conjecture, in which Oesterl\'e was interested; our goal is to say just enough to convey a decent idea of what the ``minimal discriminant'' and ``conductor'' of an elliptic curve are. The reader can, if desired, skip the next four paragraphs and jump straight to the punch line.

A general cubic plane curve is given by the equation $y^2+a_1xy+a_3y = x^3+a_2x^2+a_4x+a_6$, and sometimes by other forms, such as equation~\eqref{eqn: elliptic curve}; but we will focus on cubic curves in ``short Weierstrass form''
\begin{equation} \label{short}
y^2 = x^3+a_4x+a_6.
\end{equation}
It is always possible to find a change of variables to write a cubic plane curve in short Weierstrass form. (For example, a change of variables transforms the equation~\eqref{eqn: elliptic curve} into the form $y^2 = x^3 - 432d^2$.) In this situation, the {\em discriminant} of the curve is the quantity $\Delta = -16(4a_4^3+27a_6^2)$. If $\Delta=0$, then the cubic curve has a singularity, which is typically a node (where the graph of the curve crosses itself) but is a cusp if $a_4=a_6=0$, when the equation is simply $y^2=x^3$. But as long as $\Delta\ne0$, the cubic curve has no singularities and is called an {\em elliptic curve}.

When $a_4$ and $a_6$ are rational numbers, a change of variables can be uniquely chosen so that the coefficients $a_4$ and $a_6$ become integers with $a_4$ not divisible by the fourth power of any prime; the resulting equation is a {\em minimal model}, and its discriminant the {\em minimal discriminant}, for the elliptic curve. This minimal discriminant $\Delta$ is equal to the original discriminant times the twelfth power of a rational number, chosen so that the resulting product is an integer not divisible by the twelfth power of any prime.

Once we have a minimal model for an elliptic curve over the rational numbers, we can {\em reduce} the elliptic curve modulo any prime $p$: we simply consider the constants and variables in the equation $y^2=x^3+a_4x+a_6$ to be elements of $\Z/p\Z$, the finite field with $p$ elements. The minimal discriminant $\Delta$ over this finite field is simply the residue class of the integer $\Delta$ modulo~$p$; in particular, the reduction-at-$p$ of the elliptic curve is nonsingular (hence still an elliptic curve) precisely when $p$ does not divide~$\Delta$ (we say that the curve has {\em good reduction} at~$p$). Whenever $p$ divides $\Delta$, we say that the elliptic curve has {\em bad reduction} at~$p$. While it makes no geometric sense to talk about nodes or cusps of the ``graph'' of the elliptic curve modulo $p$---there are just a finite number of possible points, not a whole continuum---we can still categorize possible singularities algebraically, as above, into two types of bad reduction: the reduction-at-$p$ has a node (which we call {\em multiplicative reduction}) when $p$ divides $\Delta$ but not $a_4a_6$, while it has a cusp (which we call {\em additive reduction}) when $p$ divides all of $\Delta$, $a_4$, and $a_6$. (We are intentionally neglecting the more complicated cases when $p=2$ and~$p=3$.)

Finally, the {\em conductor} $N$ of an elliptic curve is a number whose prime factors are precisely those modulo which the elliptic curve has bad reduction. More specifically, $N=\prod_{p} p^{f_p}$, where the product is over all primes $p$, and $f_p$ equals $0$ if the elliptic curve has good reduction at~$p$, $1$ if it has multiplicative reduction, and $2$ if it has additive reduction. Since the primes of bad reduction are precisely the primes dividing the nonzero integer $\Delta$, all but finitely many of the $f_p$ equal $0$, and so $N$ is a well-defined positive integer. Indeed, $N$ is a multiple of $R(\Delta)$, the radical of the discriminant, and also a divisor of $R(\Delta)^2$. An elliptic curve with no primes of additive reduction is called {\em semistable}; we see that semistability is equivalent to $N=R(\Delta)$. (The breakthrough by which Andrew Wiles proved Fermat's last theorem was showing that every semistable elliptic curve was associated, through $L$-functions, to a modular form in a manner specified by the ``Taniyama--Shimura conjecture'', which is now the ``Modularity theorem''.)

{\em The punch line}: In the early 1980s, L. Szpiro formulated the following conjecture relating the minimal discriminant of an elliptic curve to its conductor.

\scj \label{conj: Szpiro}
For every $\ep>0$, there exists a positive constant $S(\ep)$ such that for any elliptic curve $E$ defined by an equation with rational coefficients,
\[|\Delta| \leq S(\ep) N^{6+\ep},\]
where $\Delta$ is the minimal discriminant of $E$ and $N$ is the conductor of~$E$.
\scje

\noindent Oesterl\'{e} observed that the newly formulated $abc$ conjecture is stronger than Szpiro's conjecture: one can deduce Szpiro's conjecture from the $abc$ conjecture, but knowing Szpiro's conjecture for all $\ep>0$, one can deduce the $abc$ conjecture only when the $\ep$ in Version~\ref{constant R^1+ep} is greater than~$\frac15$ (see \cite[Chapter VIII, exercise 8.20]{Silv} and \cite[Chapter 5, Appendix ABC]{Voj}).

In fact, Oesterl\'{e} demonstrated~\cite[pages 169--170]{Ost} that the $abc$ conjecture is actually equivalent to the following modification of the Szpiro conjecture: 

\mscj \label{conj: Szpiro, modified}
For every $\ep>0$, there exists a positive constant $S'(\ep)$ such that for any elliptic curve $E$ whose minimal model is $y^2=x^3+a_4x+a_6$,
\[
\max\{|a_4|^3,a_6^2\} \le S'(\ep) N^{6+\ep},
\]
where $N$ is the conductor of~$E$.
\mscje

\noindent Since $\Delta = -16(4a_4^3+27a_6^2)$, the modified Szpiro conjecture is clearly stronger than the original; indeed, one can take $S(\ep) = 16(4+27)S'(\ep)$ and prove the original conjecture from the modified one. But it is possible, in theory, for $\Delta$ to be small only because of extreme cancellation when the hypothetically enormous numbers $4a_4^3$ and $27a_6^2$ are added together (note that $a_4$ can be negative).

The modified Szpiro conjecture is usually stated in terms of two particular invariants $c_4$ and $c_6$ of an elliptic curve, rather than the coefficients $a_4$ and $a_6$ we have used from the short Weierstrass form~\eqref{short}; these invariants and the conductor $N$ can be associated with any elliptic curve, no matter what equation originally defines it. The invariants $c_4$ and $c_6$ are special in the sense that they suffice to determine any elliptic curve $E$ up to isomorphism (indeed, $E$ can be defined by the equation $y^2=x^3-27c_4x-54c_6$). There are algorithms for computing these invariants from various defining equations. For example (see~\cite[Section 3.2]{Crem}), an algorithm of Laska--Kraus--Connell takes $c_4$ and $c_6$, computed from any model of $E$, and outputs the minimal model of~$E$; while Tate's algorithm computes, among other things, the conductor of~$E$.

\subsection{Heuristic based on a probabilistic model} \label{prof heur sec}

While relating the $abc$ conjecture to other mathematical statements is valuable, we might be comforted by having a more instrinsic reason to believe in its truth. One way that analytic number theorists hone their beliefs about how the integers work is by creating a random-variable situation that seems to model the integer phenomenon. By rigorously showing that something happens with probability $1$ in the random model, we can gain some confidence that the analogous statement really is true in the integers. The ``Cram\'er model'' of the distribution of primes (see for example~\cite{Sound}) is probably the most well-known example of this paradigm.

Here we give a probabilistic argument, adapted from Tao~\cite{Tao}, for why one should expect the $abc$ conjecture to hold. Broadly speaking, the argument asserts that if the radicals of three integers are too small, then the ``probability'' that two of the integers sum to the third is vanishingly small; this is assuming that these numbers with small radicals are ``distributed randomly''. More specifically, we will argue that the following equivalent version of the $abc$ conjecture should be true:

\abcj \label{tao version}
Suppose $\alpha,\beta,\gamma$ are positive real numbers satisfying $\alpha+\beta+\gamma < 1$. When $M$ is sufficiently large (in terms of $\alpha,\beta,\gamma$), there are no solutions to the equation $a+b=c$ with $c\ge M$, where $a,b,c$ are relatively prime positive integers satisfying $R(a) \leq M^\alpha, R(b) \leq M^\beta, R(c) \leq M^\gamma$.
\abcje

We see rather easily that Version~\ref{finitely many c^1-ep} of the $abc$ conjecture implies this new version: given $M$ and $\alpha,\beta,\gamma$ with $\alpha+\beta+\gamma < 1$, we must have
\[
R(abc) = R(a)R(b)R(c) \le M^\alpha M^\beta M^\gamma = M^{\alpha+\beta+\gamma} < c^{\alpha+\beta+\gamma},
\]
for any triple satisfying the hypotheses of Version~\ref{tao version}. But by Version~\ref{finitely many c^1-ep} with $\ep=1-(\alpha+\beta+\gamma)>0$, there can be only finitely many such triples $(a,b,c)$; simply choose $M$ larger than the largest $c$ that occurs in any of them.

It is a little more difficult to see that Version~\ref{tao version} of the $abc$ conjecture implies Version~\ref{finitely many c^1-ep}: on the face of it, the loophole phrase ``sufficiently large (in terms of $\alpha,\beta,\gamma$)'' might allow an infinite sequence of counterexamples to Version~\ref{finitely many c^1-ep} corresponding to a sequence of distinct triples $\alpha,\beta,\gamma$. However, if there were in fact an infinite sequence of counterexamples to Version~\ref{finitely many c^1-ep} for some fixed $\ep$, then the corresponding sequence of ``vector qualities'' $(\alpha,\beta,\gamma) = \big( \frac{\log c}{\log R(a)}, \frac{\log c}{\log R(b)}, \frac{\log c}{\log R(c)} \big)$ would all lie in a compact region of $\R^3$ (namely the simplex in the positive orthant defined by $x+y+z\le 1-\ep$), and hence some subsequence of them would converge to a fixed point $(\alpha_0,\beta_0,\gamma_0)$ satisfying $\alpha_0+\beta_0+\gamma_0\le1-\ep$. Consequently, Version~\ref{tao version} of the $abc$ conjecture could be applied to the slightly larger point $(\alpha,\beta,\gamma)=\big((1+\ep)\alpha_0,(1+\ep)\beta_0,(1+\ep)\gamma_0\big)$ to derive a contradiction.

Now let's use a probabilistic model to probe Version~\ref{tao version} itself. We will need the following lemma, standard in analytic number theory (see~\cite{Tao}), saying that the numbers with a given radical are quite sparse:

\begin{proposition}
For every $\ep>0$, there exists a constant $T(\ep)$ such that for every positive squarefree integer $r$ and every $M>0$, there are at most $T(\ep)M^\ep$ integers less than $2M$ whose radical equals~$r$.
\end{proposition}

\noindent (This upper bound for the number of such integers will feature again, in a more precise form, in a lower bound given in Section~\ref{ST sec}.) Using this lemma we can estimate, given $\alpha,\beta,\gamma$, the number of triples of integers $(a,b,c)$ there are with $1\le a,b\le 2M$, $M<c\le 2M$ and $R(a) \leq M^\alpha, R(b) \leq M^\beta, R(c) \leq M^\gamma$. (Note that for the moment we are not paying attention to whether two of these numbers sum to the third number.) There are at most $M^\alpha \cdot M^\beta \cdot M^\gamma$ possibilities for the three radicals, even if we forget that they have to be squarefree and relatively prime. Choosing $\ep=\frac14(1-(\alpha+\beta+\gamma))>0$, the above proposition tells us that there are $T(\ep)M^\ep$ possibilities for $a$ for any given $R(a)$, and similarly for $b$ and $c$. Therefore the number of such triples is at most
\[
M^\alpha \cdot M^\beta \cdot M^\gamma \cdot \big( T(\ep)M^\ep \big)^3 = T(\ep)^3 M^{\alpha+\beta+\gamma+3\ep} = T(\ep)^3 M^{1-\ep}.
\]

Now, instead of looking at the specific collection of triples described in the previous paragraph, let us suppose that {\em we choose the same number of triples completely at random} from the set of triples $(a,b,c)$ with $1\le a,b\le 2M$, $M<c\le 2M$; what is the probability that at least one of the chosen triples satisfies $a+b=c$? The probability of a single chosen triple satisfying $a+b=c$ is at most $1/M$, since there is at most one correct choice out of $M$ for $c$, no matter what $a$ and $b$ are. Therefore the probability that we obtain such a chosen triple is at most $T(\ep)^3 M^{1-\ep} \cdot M^{-1} = T(\ep)^3 M^{-\ep}$. Because we have no reason to think that the actual triples described in the previous paragraph are any more or less likely to satisfy $a+b=c$ than randomly chosen triples, we are persuaded of the following heuristic: the ``probability'' is at most $T(\ep)^3 M^{-\ep}$ that there exists a triple $(a,b,c)$ with $1\le a,b\le 2M$, $M<c\le 2M$ that successfully satisfies the conditions in the second sentence of Version~\ref{tao version} of the $abc$ conjecture.

Every triple with $a+b=c$ satisfies $1\le a,b\le 2^{k+1}$, $2^k<c\le 2^{k+1}$ for a unique positive integer $k$. Let the ``$k$th event'' be the assertion that there exists a successful triple for $M=2^k$ as described above. The ``probability'' of the $k$th event, by the above heuristic, is at most $T(\ep)^3 (2^k)^{-\ep}$. Notice that the series
\[
\sum_{k=0}^\infty T(\ep)^3 (2^k)^{-\ep} = \frac{T(\ep)^3}{1-2^{-\ep}}
\]
converges to a finite number. Therefore, by the Borel--Cantelli lemma \cite[pages 51--52]{Klk}, with probability $1$ only finitely many of the events occur. We conclude that heuristically, only finitely many triples should successfully satisfy the conditions in the second sentence of Version~\ref{tao version} of the $abc$ conjecture; so if we choose $M$ large enough, we believe that there will be no counterexamples remaining.

In Section~\ref{refinements sec}, we discuss a refinement of this heuristic that enabled the authors of \cite{RST} to propose a stronger, even more precise version of the $abc$ conjecture.

\medskip

It is tempting to think of the $abc$ conjecture as a vast conspiracy that arranges the numbers with small radicals so precisely that no two of them ever add to a third. This temptation is even stronger when we see examples like the infinite families from Section~\ref{sec: inf families}: we forced the smallest and largest numbers to have tiny radicals, but the $abc$ conjecture seems, like some mystical force, to keep the middle number from being divisible by too large a perfect square. However, the above heuristic (with, say, $\alpha=\gamma=\ep$ and $\beta=1-3\ep$) offers an explanation: these families do not contain all that many triples, and probability simply dictates that it is overwhelmingly unlikely for any of the middle numbers in such a sparse set to have a small radical.

We conclude this section by remarking that the same heuristic suggests, when $\alpha+\beta+\gamma>1$, that there do exist infinitely many triples $(a,b,c)$ with $a+b=c$ and $c\ge M$ that satisfy $R(a)\le M^\alpha$, $R(b)\le M^\beta$, and $R(c)\le M^\gamma$. In other words, it is not just the numerical quality $q(a,b,c) = \frac{\log c}{\log R(abc)}$ that flirts with the boundary $q(a,b,c)=1$: the ``vector quality'' $\big( \frac{\log c}{\log R(a)}, \frac{\log c}{\log R(b)}, \frac{\log c}{\log R(c)} \big)$ actually flirts with the triangular boundary $T = \{ x,y,z\ge 0\colon x+y+z=1\}$ in~$\R^3$ at every single point. This observation leads to another question: can one find a construction of an infinite sequence of $abc$ triples such that their ``vector qualities'' approach, in the limit, a point of $T$ other than $(0,1,0)$ or one of the other two corners?

\section{Generalizations, refinements, and the state of the art} \label{sec: generalizations, results}

In this last section, we explore our most state-of-the-art knowledge about the $abc$ conjecture. We describe some (less elementary) constructions of $abc$ triples with much better quality than the ones we have seen so far; we present several refinements of the conjecture, which attempt to decide where in the space between the ``simplistic'' and actual $abc$ conjectures the exact boundary lies; we generalize the $abc$ conjecture to other settings and to more variables; and finally we discuss how close we are to actually proving the inequalities that the $abc$ conjecture asserts.

\subsection{Best known $\bm{abc}$ triples} \label{ST sec}

We already saw in Section~\ref{sec: inf families} that having examples of families of $abc$ triples helps us probe how sharp (indeed, how true) the $abc$ conjecture is. Presumably, thinking more deeply about how to find good $abc$ triples would yield even smaller radicals than the ones we have seen thus far. Stewart and Tijdeman~\cite[Theorem 2]{StwTdm} did exactly this in 1986: they came up with a construction of infinitely many $abc$ triples of higher quality than those in Section~\ref{sec: inf families}.

They proved that for any $\delta > 0$, there exist infinitely many triples $(a,b,c)$ of relatively prime positive integers with $a+b=c$ that satisfy
\begin{equation} \label{thm: StwTdm2}
c>R(abc)\exp \bigg((4-\delta)\frac{\sqrt{\log R(abc)}}{\log \log R(abc)}\bigg).
\end{equation}
All of these exps and logs can be a bit daunting for those not used to this game. With a little cunning, one may show that for any $B>1$, every sufficiently large $abc$ triple satisfying the bound~\eqref{thm: StwTdm2} also satisfies $R(abc) < c/(\log c)^B$. In fact, \eqref{thm: StwTdm2} is equivalent to a lower bound for the quality of the form
\begin{equation} \label{ST quality}
q(a,b,c) > 1 + \frac{4-\delta}{\sqrt{\log c} \cdot \log\log c}.
\end{equation}
Both of these observations show that these $abc$ triples are far better than the ones in Section~\ref{sec: inf families}. Certainly they amply refute the ``simplistic $abc$ conjecture'' and show that the epsilons in the statement of the real $abc$ conjecture must be present. But again, the lower bound for the quality tends to $1$ as $c$ grows large, and so even this better construction does not disprove the actual conjecture.

Stewart and Tijdeman's construction is quite illuminating. First, given a positive integer $r$ and a parameter $X$ that is far larger than $r$, they consider the set of integers up to $X$ whose factorizations contain only the first $r$ odd primes $p_1,\dots,p_r$; such integers are called ``$p^r$-friable'' (or ``$p^r$-smooth''). They obtain good estimates for the number of such integers by noting that the number of solutions to $p_1^{n_1}\cdots p_r^{n_r} \le X$ is the same as the number of lattice points $(n_1,\dots,n_r)$ in the $r$-dimensional simplex (high-dimensional pyramid)
\[
\big\{ x_1,\dots,x_r\ge0 \colon x_1\log p_1 + \cdots + x_r\log p_r \le \log X \big\}.
\]
When $X$ is large, this number of lattice points is essentially the $r$-dimensional volume of the simplex, which is easy to calculate. Their resulting lower bound for the number of $p^r$-friable integers up to $X$ is a more precise version of the Proposition from Section~\ref{prof heur sec}.

Then, they point out that two of these $p^r$-friable integers must be congruent to each other modulo a high power of $2$; indeed, if the number of such integers exceeds $2^k$, then the pigeonhole principle forces two of them to lie in the same residue class modulo~$2^k$. (And if those two happen to have common factors, their quotients by their greatest common divisor will also be congruent to each other.) These two integers and their difference form a triple satisfying $a+b=c$; the product of the radicals of the first two integers is at most $p_1\cdots p_r$; and the radical of their difference is at most $X/2^{k-1}$.

Of course, it is necessary to write down explicitly the relationships among all of these functions and parameters; analytic number theorists learn tools that are precisely suited for converting from the above sketch to a full quantitative proof. When the dust settles, the inequality~\eqref{thm: StwTdm2} is the payoff.

Later, van Frankenhuysen~\cite{vF} added an improvement to the argument of Stewart and Tijdeman: instead of using the full lattice of integer points $(n_1,\dots,n_r)$, he chose a sublattice sitting askew inside the full integer lattice in such a way that the points in the sublattice were relatively more tightly packed together, in the same way that a pyramidal stack of oranges in the grocery store takes up less space than they would if we insisted on placing each one directly atop the one below it. In this way, and using high-dimensional sphere-packing bounds already in the literature, he showed that one could improve the constant $4-\delta$ in the above inequalities to $6.068$.

At the end of Section~\ref{sec: backgr, motiv}, we described the vector quality of an $abc$ triple; we remark here that even these fancy $abc$ triples of Stewart/Tijdeman and van Frankenhuysen have the property that their vector qualities converge to $(0,1,0)$ or one of the other two corners, rather than some intermediate point $(\alpha,\beta,\gamma)$ with $\alpha+\beta+\gamma=1$ and $\alpha\beta\gamma>0$. In other words, we still don't know how to construct an infinite family of $abc$ triples where each of the three integers has a radical that is significantly smaller than itself.

\subsection{Refinements of the $\bm{abc}$ conjecture} \label{refinements sec}

In a sense, the loophole phrases ``only finitely many'' and ``there exists a positive constant'' make it hard to actually determine from data whether the $abc$ conjecture is acting like the truth. In 1996, Baker~\cite{Bkr} refined the conjecture to provide some insight into how the constant $K(\ep)$ in Version~\ref{constant R^1+ep} of the $abc$ conjecture should depend on~$\ep$. Let the function $\omega(n)$ denote the number of distinct primes (that is, ignoring repetitions) dividing~$n$. Baker proposed the following refinement: there exists an absolute constant $K_1$ such that, if $a,b,c$ are relatively prime positive integers satisfying $a+b=c$, then 
\[
c < \ep^{-K_1\min\{\omega(ab),\omega(ac),\omega(bc)\}}R(abc)^{1+\ep}
\]
for any $\ep>0$. (The minimum in the exponent is there to help us: we can keep whichever two of the three numbers have the fewest prime factors between them.) Although this bound has an extra dependence on $a,b,c$, this dependence turns out to be smaller than $R(abc)^\ep$, and so we have not significantly altered the shape of the conjecture. It is true that there is still an unknown absolute constant $K_1>0$ in this formulation; but at least now this constant is independent of $\ep$, which is often quite important when making deductions.

Baker did demonstrate that there exist infinitely many $abc$ triples satisfying the related inequality
\[
c > K_2 \ep^{1-\min\{\omega(ab),\omega(ac),\omega(bc)\}}R(abc)^{1+\ep}
\]
for some absolute constant $K_2>0$. In fact, his proof relies upon estimates for ``linear forms in logarithms'', a profound technical tool in Diophantine approximation for which Baker was awarded the Fields Medal in 1970. He also mentions that Granville conjectured that there exists an absolute constant $K_3>0$ such that $c < K_3^{\Omega(abc)} R(abc)$ for all $abc$ triples, where $\Omega(n)$ counts the number of prime factors of $n$ with multiplicity (so for example, $\omega(72)=2$ but $\Omega(72)=5$). Notice that there is no exponent $1+\ep$ on the right-hand side of this conjecture! But of course the slack has to go somewhere: $\Omega(abc)$ can be far larger than $\omega(abc)$.

Most recently, another refinement has been put forward by Robert, Stewart, and Tenenbaum~\cite{RST}, following up on a heuristic proposed by van Frankenhuysen in his PhD thesis. In Section~\ref{prof heur sec} we described the heuristic assumption that statistically, $R(c)$ is distributed independently from $R(a)$ and $R(b)$ when $a, b, c$ are relatively prime and $a+b=c$; we then estimated how many integers up to $M$ have radicals bounded by $M^\alpha$ and so on. Through an extremely careful study of the function that counts how many integers up to $M$ have their radical bounded by a second parameter $Y$, Robert, Stewart, and Tenenbaum proposed the following two-part conjecture: First, there exists a real number $K_4$ such that all $abc$ triples satisfy
\begin{equation} \label{conj: RST}
c < R(abc)\exp\bigg(4\sqrt{\frac{3\log R(abc)}{\log \log R(abc)}}\bigg(1+\frac{\log \log \log R(abc)}{2 \log \log R(abc)}+\frac{K_4}{\log \log R(abc)}\bigg)\bigg);
\end{equation}
second, there exists a real number $K_5$ such that infinitely many $abc$ triples satisfy
\begin{equation} \label{conj: RST2}
c> R(abc)\exp\bigg(4\sqrt{\frac{3 \log R(abc)}{\log \log R(abc)}}\bigg(1+\frac{\log \log \log R(abc)}{2\log \log R(abc)}+\frac{K_5}{\log \log R(abc)}\bigg)\bigg).
\end{equation}
(To tell the truth, they included more detailed versions of this refinement with even more logs in the picture!) Note that the right-hand side of their second conjecture~\eqref{conj: RST2} is a little bit larger than the lower bound~\eqref{thm: StwTdm2} coming from the construction of Stewart and Tijdeman, since we are dividing only by $\sqrt{\log\log R(abc)}$ in the main term inside the exponential instead of $\log\log R(abc)$. It is actually easier than it might seem to show that their first conjecture~\eqref{conj: RST} really does imply Version~\ref{finitely many R^1+ep} of the $abc$ conjecture.

Interestingly, a special quantity arises from these conjectures of Robert, Stewart, and Tenenbaum. Define the ``merit" of an $abc$ triple to be
\[
m(a,b,c) = \big(q(a,b,c)-1\big)^2\log{R(abc)}\log{\log{R(abc)}}.
\]
Every infinite family of $abc$ triples ever established has the property that their merit tends to~$0$ (it is not hard to verify this for the Stewart--Tijdeman examples, for instance, assuming that the right-hand side of the inequality~\eqref{ST quality} is in fact the correct size of the quality). If the $abc$ conjecture were false, it is an easy deduction from Version~\ref{quality version} that the merit would be unbounded above. But it would actually follow from the conjectures~\eqref{conj: RST} and~\eqref{conj: RST2} that the lim sup of all merits of all $abc$ triples equals $48$ exactly! So the merit is somehow an incredibly fine-scale measurement of an $abc$ triple, one that looks at the boundary between possible and impossible through a powerful microscope.

For the record, the largest merit found to date is approximatey $38.67$, which comes from the $abc$ triple $(2543^4 \cdot 182587 \cdot 2802983 \cdot 85813163, 2^{15} \cdot 3^{77} \cdot 11 \cdot 173, 5^{56} \cdot 245983)$ discovered by Ralf Bonse in 2011. de Smit's web site~\cite{Smt} lists the $131$ known $abc$ triples with merit greater than~$24$.

\subsection{Other alterations of the conjecture} \label{other alterations sec}

In addition to these results, some interesting generalizations and refinements of the $abc$ conjecture have appeared in various contexts.

\subsubsection*{Congruence $abc$ conjecture}

First we state yet another version of the conjecture, which looks like it concerns only a small subset of $abc$ triples but is actually equivalent to the other versions we have seen so far.

\abcj \label{congruence}
For every positive integer $N$ and every $\ep>0$, there exists a positive constant $E(N,\ep)$ such that all triples $(a,b,c)$ of relatively prime positive integers with $a+b=c$ and $N\mid abc$ satisfy \[c \le E(N,\ep) R(abc)^{1+\ep}.\]
\abcje

The special case $N=1$ is of course the familiar $abc$ conjecture (specifically, Version~\ref{constant R^1+ep}). Motivated by Oesterl\'{e}'s observation \cite{Ost} that the special case $N=16$ of this new version implies the full $abc$ conjecture, Ellenberg \cite{Ell} demonstrated that in fact $16$ can be replaced by any integer $N$, thus showing that Version~\ref{congruence} of the $abc$ conjecture really is equivalent to the others.

\subsubsection*{Uniform $abc$ conjecture for number fields} \label{number field sec}
An important focus in algebraic number theory is the study of {\em number fields}, which are finite field extensions of the field $\Q$ of rational numbers. Equivalently, a number field is a field of the form $\Q(\alpha)$, where $\alpha$ is a root of a polynomial with integer coefficients (an {\em algebraic number}). For example, $\alpha=i+\sqrt2$ is an algebraic number since it is a root of the polynomial $x^4-2 x^2+9$; consequently, $\Q(\alpha)$ is a number field, consisting of all complex numbers of the form $r+s\alpha+t\alpha^2+u\alpha^3$ for rational numbers $r,s,t,u$. Arithmetic can be done in a consistent way in number fields, almost as nicely as in the rational numbers themselves, and their study is essential to our understanding of solutions of polynomial equations.

Vojta \cite[page 84]{Voj} formulated a generalization of the $abc$ conjecture to number fields, pointing out many notable consequences of this generalization (see also \cite{Bom, Elk, GrvStk}). However, that formulation contains some unfamiliar terminology that would be too laborious to define here. The next paragraph, therefore, is intended for those who are more familiar with algebraic number theory; other readers may skip that paragraph and at least get an impressionistic idea of the statement of the uniform $abc$ conjecture.

Let $K/\Q$ be a number field of degree $n$ with discriminant $D_K$. For each prime ideal $\mathfrak{p}$ of $K$, let $|\;|_{\mathfrak{p}}$ be the corresponding $\mathfrak{p}$-adic absolute value, normalized so that $|\mathfrak{p}|_{\mathfrak{p}} = \Nm_{K/\Q}(\mathfrak{p})^{-1/n}$; for each real or complex embedding $\tau$ of $K$, let $|\alpha|_\tau = |\tau(\alpha)|^{1/n}$ be the corresponding normalized Archimedean absolute value, where $|\;|$ is the modulus of a complex number. Let the height of the $m$-tuple $(\alpha_1, ..., \alpha_m) \in K^m$ be given by
\[
H(\alpha_1, ..., \alpha_m) = \prod_{v} \max\big\{|\alpha_1|_v, ..., |\alpha_m|_v\big\},
\]
where the product goes over all places $v$ (prime ideals and embeddings). Finally, let the conductor of the $m$-tuple be given by
\[
N(\alpha_1, ..., \alpha_m) = \prod_{\mathfrak{p} \in I}^{} |\mathfrak{p}|_{\mathfrak{p}}^{-1},
\]
where $I$ is the set of prime ideals $\mathfrak{p}$ such that $|\alpha_1|_{\mathfrak{p}}, \dots, |\alpha_m|_{\mathfrak{p}}$ are not all equal. Then we have:

\ucj
For every $\ep>0$, there exists a constant $U(\ep)>0$ with the following property: for every number field $K$ of degree $n$ over $\Q$, and every triple $(a,b,c)$ of elements of $K$ satisfying $a+b+c=0$,
\[
H(a,b,c) \leq U(\ep) \big(D_K^{1/n}  N(a,b,c)\big)^{1+\ep}.
\]
\ucje

To shed some light on the relationship between this number field version and the usual $abc$ conjecture, we remark that if $K=\Q$, then $n=1$ and $D_K=1$; furthermore, if $(a,b,c)$ is a relatively prime triple of integers, then the height $H(a,b,c)$ is simply $\max\{|a|,|b|,|c|\}$ and the conductor $N(a,b,c)$ is simply $R(|abc|)$. (The above definitions of the height and conductor have the convenient property that they do not change if every element of the $m$-tuple is multiplied by the same factor, and so a relative primality hypothesis is actually unnecessary for this generalization.) Therefore the $K=\Q$ case of the uniform $abc$ conjecture is exactly Version~\ref{constant R^1+ep} of the $abc$ conjecture, once we take absolute values of the three numbers and reorder them so that $c$ is the largest.

\subsubsection*{Additional integers}

Our last generalization, stated by Browkin and Brzezinski~\cite{BB}, incorporates more variables than the three we have been working with so far.

\ncj
For every integer $n\ge3$ and every $\ep>0$, there exists a positive constant $B(n,\ep)$ such that all relatively prime $n$-tuples $(a_1,\dots,a_n)$ of nonzero integers with $a_1+\cdots+a_n=0$ and no vanishing subsums satisfy \[\max\{|a_1|,\dots,|a_n|\} \le B(n,\ep) R(|a_1\cdots a_n|)^{2n-5+\ep}.\]
\ncje

\noindent Here, ``no vanishing subsums'' means that it is not possible to reorder $a_1,\dots,a_n$ so that $a_1+\dots+a_k=0=a_{k+1}+\cdots+a_n$ for some $1\le k\le n-1$; this hypothesis is necessary because of trivial examples such as $(a_1,a_2,a_3,a_4) = (2^n,-2^n,3^n,-3^n)$, which is a relatively prime quadruple even though some pairs of terms have huge common factors. This $n$-variable version is our familiar friend when $n=3$: given such a triple $a_1,a_2,a_3$, one of them has a different sign than the other two, and we recover Version~\ref{constant R^1+ep} of the $abc$ conjecture by letting $c$ be the absolute value of the one with a different sign and $a,b$ the absolute values of the other two.

Browkin and Brzezinski constructed examples showing that the exponent $2n-5+\ep$ on the right-hand side cannot be reduced; their constructions are rather similar to the transfer method described in Section~\ref{transfer section}. Taking $n=4$ for example, if $(a,b,c)$ is any $abc$ triple, we may set
\begin{equation} \label{BB example 4}
(a_1,a_2,a_3,a_4) = (a^3,3abc,b^3,-c^3),
\end{equation}
which one can check does satisfy $a_1+a_2+a_3+a_4=0$; for these quadruples,
\[
\max\{|a_1|,\dots,|a_n|\} = c^3 \ge \big(R(abc)\big)^3 \ge \big(\tfrac13R(|a_1a_2a_3a_4|)\big)^3 = \tfrac1{27}\big(R(|a_1a_2a_3a_4|)\big)^{2\cdot4-5}.
\]
Equation~\eqref{BB example 4} is the $n=4$ case of a sequence of impressive identities: when $n\ge3$,
\begin{equation} \label{BB examples}
a_j = \frac{2n-5}{2j-1} \binom{n+j-4}{2j-2} a^{2j-1} (bc)^{n-j-2} \text{ for } 1\le j\le n-2, \quad a_{n-1} = b^{2n-5}, \quad a_n = -c^{2n-5}
\end{equation}
is an $n$-tuple satisfying the hypotheses of the $n$-variable version of the $abc$ conjecture. For these $n$-tuples, the maximum absolute value is $c^{2n-5}$, while $R(|a_1\cdots a_n|)$ is at most a constant (the product of all the primes up to $2n-5$, say) times $R(abc)$, which is at most $c$ times a constant when $(a,b,c)$ is an $abc$ triple. Not only does this construction show that the exponent $2n-5+\ep$ would be best possible, it also shows that the $n$-variable version for any $n\ge4$ implies the usual three-variable $abc$ conjecture. (It seems less clear whether, for example, the $5$-variable version of the $abc$ conjecture implies the $4$-variable version.)

Interestingly, the statement of the $n$-variable version of the $abc$ conjecture is not what one would predict from a probabilistic heuristic like the one described in Section~\ref{prof heur sec}: the analogous argument would lead again to a conjecture with exponent $1+\ep$ on the right-hand side. In this case, probability would lead us astray---but presumably because the set of counterexamples is extremely thin, coming only from constructions like equation~\eqref{BB example 4}. In fact, it follows from a sufficiently strong version of Vojta's conjecture \cite[Conjecture 2.3]{Voj2} that the exponent $2n-5+\ep$ can be reduced to $1+\ep$ if we exclude a finite number of constructions like equation~\eqref{BB examples} for each $n$. Even without excluding these constructions, it might be possible to reduce the exponent somewhat if we insist that the $n$-tuple be pairwise relatively prime, rather than just relatively prime as a tuple. In fact, in the function field case \cite{BM} such theorems have been worked out under intermediate assumptions such as every $m$-subtuple of the $n$-tuple being relatively prime; these theorems could serve as motivation for analogous versions of the $n$-variable conjecture, in the spirit of the previous section.

\subsection{Progress towards the $\bm{abc}$ conjecture} \label{sec: current status}

The first players of this game were Stewart and Tijdeman \cite{StwTdm}: in 1986 they proved that
\begin{equation}
\label{thm: StwTdm1}
c < \exp\big(K_6 R(abc)^{15}\big)
\end{equation}
for some constant $K_6>0$. Their proof used bounds on linear forms in logarithms similar to those mentioned in Section~\ref{refinements sec}, in particular a $p$-adic version due to van der Poorten. Subsequently, Stewart and Yu \cite{StwYu1,StwYu2} improved the bound~\eqref{thm: StwTdm1} to
$$
c < \exp\big(K_7 R(abc)^{{1/3}}(\log R(abc))^3 \big)
$$
for some constant~$K_7>0$. They achieved this improvement by replacing van der Poorten's $p$-adic estimates with even stronger ones due to Yu. Despite the fact that these results were hard-earned and at least bring the problem of bounding $c$ into the realm of the finite, neither inequality is as good as $c < R(abc)^B$ for any fixed~$B$.

The number theory community has been abuzz with the topic of the $abc$ conjecture the past few years. In August 2012, Shinichi Mochizuki released the final installment of his series of four papers on ``inter-universal Teichm\"{u}ller theory'', in which he claimed to have proven the $abc$ conjecture as a consequence of his work. His proof, with its incredible length and heavy dependence on his past work in anabelian geometry---a new and untested field with a limited number of practitioners---is still under verification by the mathematical community. Moreover, due to the introduction of several arcane objects such as ``Frobenioids," ``log-theta-lattice", and ``alien arithmetic holomorphic structures," a cautious response from the mathematical community was inevitable.

Mochizuki published a progress report in December 2013, informing the community of the advencement that had been made towards verifying his results. (See the Polymath page~\cite{poly} for useful links to Mochizuki's papers, progress report, announcements, and other related topics.) Members of his home university have studied his preparatory papers and waded through his manuscripts on inter-universal Teichm\"{u}ller theory, communicating with Mochizuki on suggested improvements and adjustments to be made; they plan to give seminars on the material starting in the fall of 2014. On the other hand, due to the esoteric nature of Mochizuki's work (and the presence of some at least superficial mistakes in the deduction of the $abc$ conjecture from his theory), it has been hard for others to attest to the validity of his results. While a wave of colleagues around the world was drawn to the task of understanding his exotic, potentially revolutionary work, the reality is that it is difficult for most academics to pause their own research to invest the necessary energy. Several skilled mathematicians spent a good deal of time trying to understand how the arguments were structured but, after making little headway in being able to verify Mochizuki's claims, eventually abandoned the project.

This unsettled state of affairs begs the question: when does one say that a problem in mathematics has been solved? Many of us would like to think we have an absolute standard, where proofs are accepted if and only if they are completely rigorous and complete, line by line, like a successfully compiling computer program. But in practice, we tolerate typos, allusions to proofs of similar cases, sketches of arguments, and occasional exercises for the reader as acceptable parts of research papers; our standard of proof in mathematics is a social construct \cite[Section 4]{Thu}. Researchers in specialized fields form their own epistemic communities and move forward in clusters, building around one another's work, and sharing their knowledge with researchers in neighboring areas as they can.

A mathematician's results, then, are accepted only when her primary audience---the cohort of experts occupying the same niche---has validated their accuracy. In this case, with Mochizuki's original and complex work, it will take some time for more mathematicians to surmount the barrier and begin exporting the ideas to the wider community. In the best possible world, experts will come to agree that the papers contain a proof to one of the most significant problems in number theory, as well as the foundations of new areas of research. But until and unless that happens, we must be content with the $abc$ conjecture remaining a mystery, at least for now.

\section*{Acknowledgments}

The authors are happy to thank Tim Dokchitser for indicating where to locate the material in Section~\ref{j section} in the literature; Carlo Beenakker and Cam Stewart for their help in finding attributions for some of the examples in Section~\ref{folklore sec}; and joro, Felipe Voloch, and Trevor Wooley for helpful comments on the $n$-variable version of the $abc$ conjecture mentioned in Section~\ref{other alterations sec} and the connection to Vojta's work.



\end{document}